\documentclass[conference]{IEEEtran}
\IEEEoverridecommandlockouts
\usepackage{cite}
\usepackage{amsmath,amssymb,amsfonts}
\usepackage{algorithmic}
\usepackage{graphicx}
\usepackage{textcomp}
\usepackage{xcolor}
\def\BibTeX{{\rm B\kern-.05em{\sc i\kern-.025em b}\kern-.08em
    T\kern-.1667em\lower.7ex\hbox{E}\kern-.125emX}}
\begin{document}

\title{A conjecture related to the Newman phase\\
%\thanks{This work was done outside of Amazon and the views are that of the author's only.}
}

\author{\IEEEauthorblockN{Shravan Mohan}
\thanks{This work is not related to the work of the author at Amazon.}
\IEEEauthorblockA{Mantri Residency, Bengaluru\\
shravan.rammohan@gmail.com}
}

\maketitle

\begin{abstract}
A conjecture is proposed concerning the recovery of a discrete magnitude spectrum through a nonlinear transformation involving the Newman's phase sequence. Given a discrete magnitude spectrum sampled from a continuous function, consider the process of applying a complex exponential with the Newman's phase sequence, computing the inverse discrete Fourier transform (IFFT), taking the absolute value of the result, and reversing the time-domain signal. The conjecture states that Newman's phase sequence, defined by a formula $\phi^{(N)}[k] = \frac{\pi (k-1)^2}{N}$,  asymptotically recovers the original magnitude spectrum as the number of samples $N \to \infty$. Notably, the phase sequence is also independent of the input signal and is unique up to an overall constant phase shift. The broader implications of this conjecture remain to be fully understood, but the phenomenon raises fundamental questions about the role of phase in nonlinear spectral recovery.
\end{abstract}


\section{Introduction}

\noindent
In the course of exploring the effect of phase on the structure of time-domain signals, an empirical phenomenon was observed. Given a non-negative continuous function defined on the interval \( [-\pi, \pi] \), various discrete magnitude spectra obtained via uniform sampling were considered. These spectra were chosen to span a wide range of behaviors—symmetric and asymmetric, smooth and non-smooth. When the corresponding phase spectrum was chosen as the \emph{Newman phase sequence}, the pointwise magnitude of the inverse discrete Fourier transform was observed to resemble the \emph{reversed} version of the input magnitude spectrum, increasingly so as the number of frequency samples \( N \) grew. Note that the Newman phase sequence is defined as follows: for each \( N \geq 1 \) and \( k \in \{0, 1, \dots, N-1\} \),
\[
\phi_{\text{Newman}}(k, N) = \pi \frac{k^2}{N}.
\]

\medskip

\noindent
To correct this spatial reversal, a time-reversal step is necessary. It is important to note that adding a linear phase term \( \pi k \) 
to the Newman phase, equivalent to time-shifting, does not correct the reversal. This is also evident with the example shown in Figure. 2. With this modification, the magnitude of the resulting time-domain signal began to closely match the original magnitude spectrum directly. In the case of non-smooth magnitude spectra—particularly those with jump discontinuities—oscillatory artifacts reminiscent of the Gibbs phenomenon were also observed. These consistent empirical observations across a wide class of input functions motivated the following conjecture.

\medskip

\begin{figure}[t]
    \centering
    \includegraphics[width=3.5in]{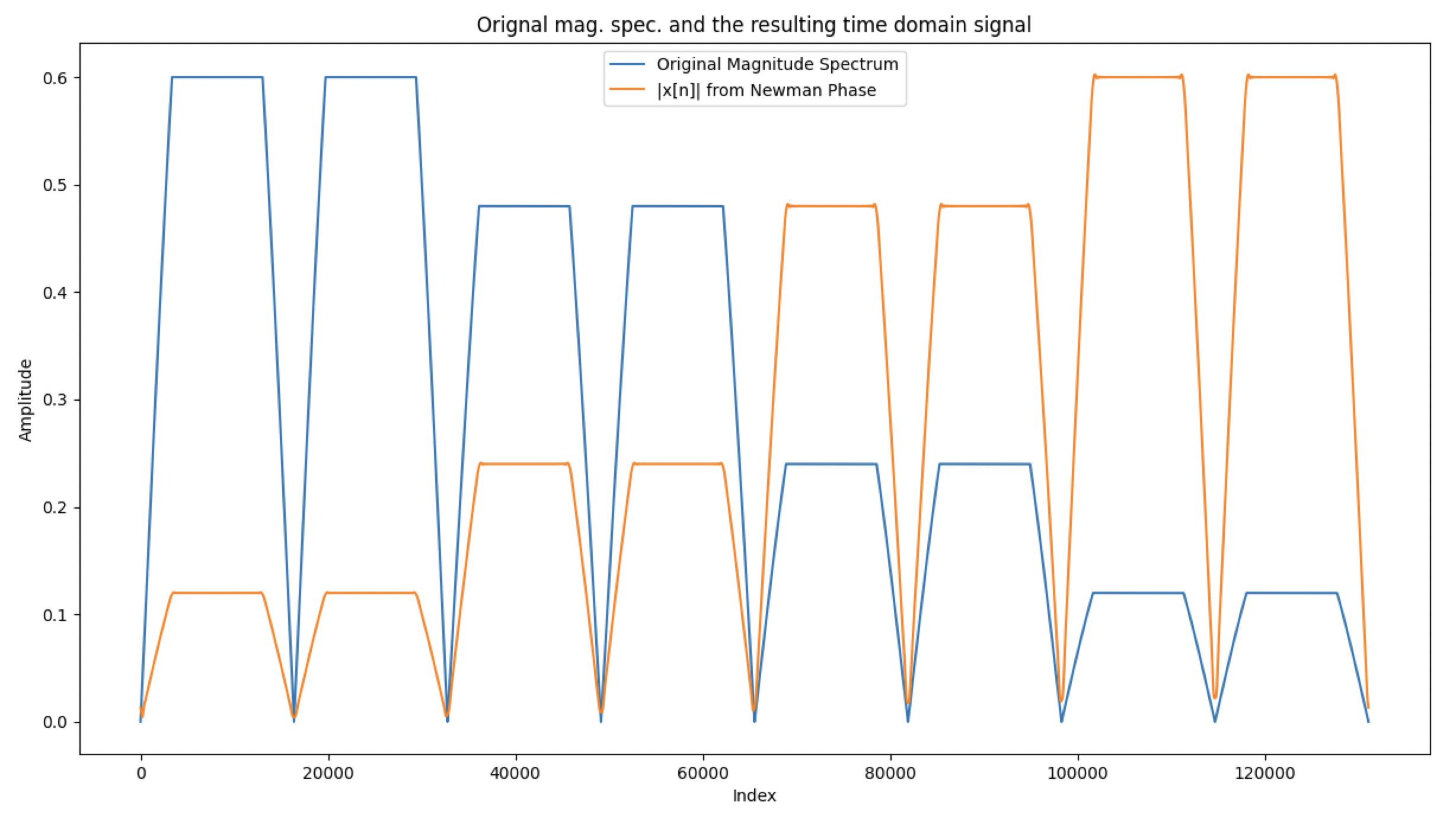}
    \caption{
        Comparison of the original magnitude spectrum (blue) and the magnitude of the time-domain signal obtained after applying the Newman phase(orange). The resulting magnitude does not align with the original spectrum, illustrating that a simple linear phase shift—equivalent to a time-domain translation—is insufficient. This motivates the need for an additional time reversal of the inverse DFT output.
    }
    \label{fig:reversal}
\end{figure}

\begin{figure*}[ht]
    \centering
    \includegraphics[width=1\linewidth]{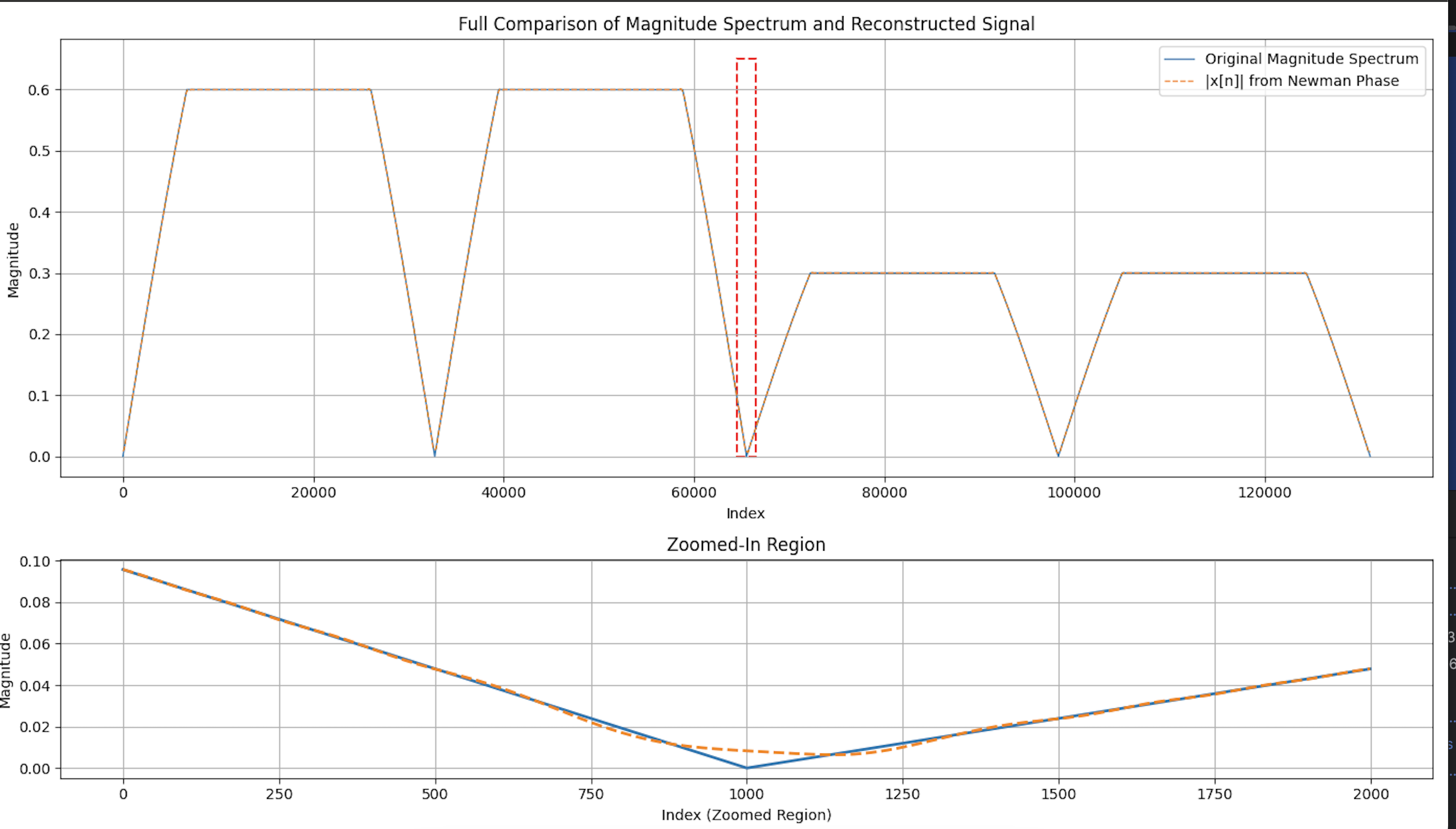}
    \caption{Comparison of the original discrete magnitude spectrum and the absolute value of the reconstructed time-domain signal obtained using the Newman phase sequence with an added linear phase term. The top plot shows the full range for \( N = 131072 \); the red dashed box highlights a region with a discontinuity. The bottom plot zooms into this region, illustrating the small but visible oscillations (Gibbs-like behavior) that appear near the discontinuity in the reconstructed signal.}
    \label{fig:zoomed_diff}
\end{figure*}

\noindent
Let \( m : [-\pi, \pi] \to \mathbb{R}_{\geq 0} \) be a non-negative continuous function. For each integer \( N \geq 1 \), define a length-\( N \) discrete magnitude spectrum by uniform sampling:
\[
M^{(N)}[k] = m\left(-\pi + \frac{2\pi k}{N}\right), \quad k = 0, 1, \dots, N-1.
\]
Let \( \Phi(k, N) = \phi_{\text{Newman}}(k, N)\), where \( \phi_{\text{Newman}}(k, N) \) denotes the Newman phase sequence. Form the complex-valued sequence:
\[
X^{(N)}[k] = M^{(N)}[k] \cdot e^{i \Phi(k, N)}.
\]
Compute the inverse discrete Fourier transform:
\[
x^{(N)}[n] = \frac{1}{N} \sum_{k=0}^{N-1} X^{(N)}[k] \cdot e^{2\pi i kn / N}, \quad n = 0, 1, \dots, N-1.
\]
Then define the pointwise magnitude:
\[
\widetilde{M}^{(N)}[k] = \left| x^{(N)}[k] \right|, \quad k = 0, 1, \dots, N-1.
\]

\medskip

\noindent
With this, the following conjecture is proposed:

\medskip

\noindent
\textbf{Conjecture.} \textit{There exists a subset \( \mathcal{D} \subset C^{\infty}([-\pi, \pi]; \mathbb{R}_{\geq 0}) \) of smooth functions such that $\mathcal{D}$ is dense in $C^{\infty}([-\pi, \pi]; \mathbb{R}_{\geq 0})$, and for every \( m \in \mathcal{D} \), the sequence \( \widetilde{M}^{(N)} \) constructed using the phase \( \Phi(k, N) = \phi_{\text{Newman}}(k, N) \) satisfies}
\[
\lim_{N \to \infty} \sup_{0 \leq k < N} \left| \widetilde{M}^{(N)}[N-k] - M^{(N)}[k] \right| = 0, ~~a.e.
\]
\textit{Moreover, the phase function \( \Phi(k, N) \) is independent of the input function \( m \), and is unique up to an additive constant phase shift.}

\section*{Genesis of the Newman phase sequence}

\noindent Let \( P(z) = \sum_{k=0}^{n} a_k z^k \) be a polynomial of degree \(n\) whose coefficients satisfy the constraint
\[
|a_k| = 1 \quad \text{for all } k = 0, 1, \dots, n.
\]
Define the quantity
\[
I_n = \min_{|a_k| = 1} \frac{1}{2\pi} \int_0^{2\pi} \left| \sum_{k=0}^{n} a_k e^{i k \theta} \right| \, d\theta,
\]
i.e., the minimum \(L^1\)-norm of such a polynomial over the unit circle. The central question is:  
\emph{What is the asymptotic behavior of \( I_n \) as \( n \to \infty \)?}\\

\noindent Earlier results using Schwarz-type inequalities had established that
\[
I_n \leq \sqrt{n+1}.
\]
Newman \cite{newman1965} aimed to determine whether the lower bound matched this up to leading order—that is, whether
\[
\lim_{n \to \infty} \frac{I_n}{\sqrt{n}} = 1.
\]
To answer this, he constructed a specific polynomial whose \(L^1\) norm is asymptotically \( \sqrt{n} - c \), thereby resolving the problem and establishing the correct order of growth. Specifically, he considers the polynomial \( P(z) = \sum_{k=0}^{n} a_k z^k \) with 
\[
a_k = \exp\left( \frac{2\pi i k^2}{n + 1} \right), \quad \text{for } k = 0, 1, \dots, n,
\]
and shows that 
\[
\frac{1}{2\pi} \int_0^{2\pi} |P(e^{i\theta})|\,d\theta \geq \sqrt{n} - c,
\]
for some absolute constant \( c \), thereby showing that the lower bound matches the known upper bound asymptotically. While the original paper does not name this construction, later literature refers to the coefficient structure \( a_k = e^{2\pi i k^2 / (n+1)} \) as the \emph{Newman phase}. This terminology honors Newman's elegant use of quadratic exponential phases to solve a difficult extremal problem in harmonic analysis. 

\begin{figure*}[t]
    \centering
    \includegraphics[width=1\linewidth]{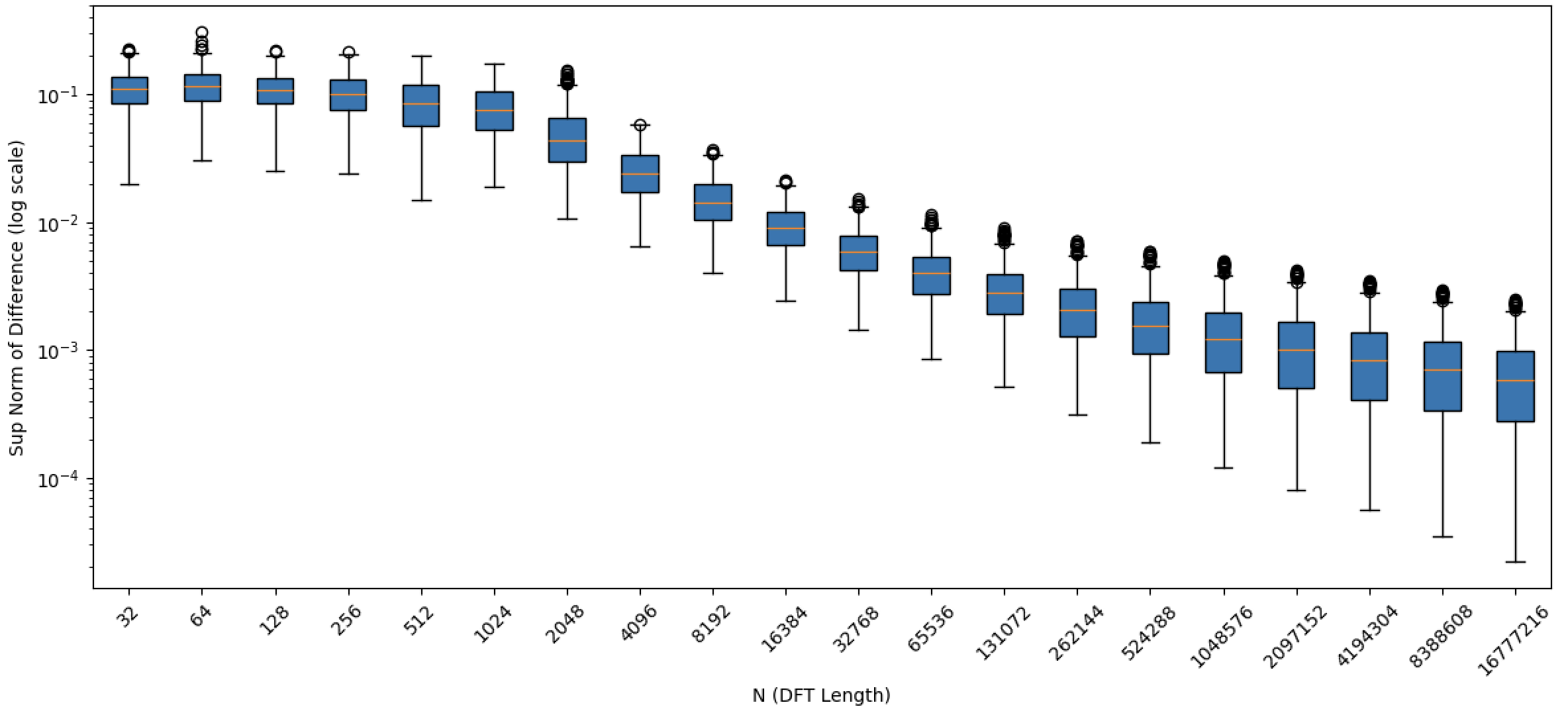}
    \caption{{Box plot of the root mean square (RMS) error 
    $\sqrt{\frac{1}{N} \sum_{k=0}^{N-1} \big( M[k] - |x[n]|[k] \big)^2}$ 
    for 1{,}000 random sum-of-sinusoids signals with random phases and a positive bias, using only the Newman phase sequence.
    For each signal, the magnitude spectrum was sampled at increasing DFT lengths $N$ (from $32$ up to $1{,}677{,}216$).
    The plot shows the distribution of reconstruction errors after applying the Newman phase, inverse discrete Fourier transform, and time reversal.
    The logarithmic scale highlights how the RMS error tends to decrease with larger $N$, supporting the empirical convergence behavior for this signal class.}
    }   
    \label{fig:zoomed_diff}
\end{figure*}
\section{Applications of Newman Phases}

\noindent Although the Newman phase construction was originally introduced in the context of extremal polynomial problems, it has since found impactful applications in engineering—particularly in the design of multitone and multicarrier signals with \emph{low Peak-to-Average Power Ratio (PAPR)}. PAPR, defined as the ratio of a signal’s peak amplitude to its RMS value, is a critical measure in communication systems such as orthogonal frequency division multiplexing (OFDM), where high PAPR can introduce nonlinear distortion and reduce amplifier efficiency. In a foundational study, Stephen Boyd \cite{boyd1986} demonstrated that using the \emph{Newman phase sequence}
\[
\theta_k = \frac{\pi k^2}{N}, \quad k = 0, 1, \dots, N-1,
\]
for a multitone signal results in a crest factor below 6 dB. These quadratic-phase sequences ensure the signal has a nearly uniform amplitude distribution, minimizing the time spent near peak levels and thereby reducing PAPR.\\

\noindent Building on this idea, Jayalath et al.\ \cite{jayalath2000} proposed using Newman phase sequences in \emph{Selected Mapping (SLM)} schemes for OFDM. Their approach yielded a \emph{PAPR reduction of over 4.5 dB} for 256-carrier systems, with lower computational complexity compared to traditional Partial Transmit Sequence (PTS) methods. In parallel, Al-Imari and Hoshyar \cite{alimari2012} extended this application to \emph{Low Density Signature OFDM (LDS-OFDM)}, a multiple access scheme. By applying Newman phases to the signature sequences, they achieved a \emph{3.6 dB PAPR reduction} compared to random phasing, without altering the LDS-OFDM system architecture. In the domain of optical communications, Dai, Guo, and Huang \cite{dai2013} showed through simulation that using Newman phase sequences in \emph{optical OFDM systems} significantly reduces PAPR and leads to an \emph{optical signal-to-noise ratio (OSNR) gain of 8 dB}. These examples illustrate just a few known engineering applications of Newman phases; there may well be others that the author is not currently aware of.

\section{Numerical insights}

\noindent
Fix a positive integer \( N \). Define the diagonal matrix \( \Lambda^{(N)} \in \mathbb{C}^{N \times N} \) with diagonal entries given by the Newman phase sequence:
\[
\Lambda^{(N)} = \mathrm{diag}\left( e^{i  \frac{\pi k^2}{N}} \right)_{k=0}^{N-1}.
\]

\noindent
Given a discrete real non-negative magnitude spectrum \( M^{(N)} \in \mathbb{R}_{\geq 0}^N \), the complex spectrum is obtained by attaching the Newman phase:
\[
X^{(N)} = \Lambda^{(N)} M^{(N)}.
\]

\noindent
Let \( F^{(N)} \in \mathbb{C}^{N \times N} \) denote the inverse discrete Fourier transform (IDFT) matrix, with entries:
\[
F^{(N)}[n, k] = \frac{1}{N} e^{2\pi i nk / N}, \quad 0 \leq n, k < N.
\]

\noindent
The time-domain signal is given by:
\[
x^{(N)} = F^{(N)} \Lambda^{(N)} M^{(N)}.
\]

\noindent
To correct for the spatial reversal observed in earlier experiments, a time-reversal operator \( R \in \mathbb{R}^{N \times N} \) can be applied, defined by:
\[
R[n, k] = \delta_{n, N - 1 - k}, \quad 0 \leq n, k < N,
\]
where \( \delta_{i,j} \) is the Kronecker delta. Applying this yields the time-reversed signal:
\[
\tilde{x}^{(N)} = R x^{(N)} = R F^{(N)} \Lambda^{(N)} M^{(N)}.
\]

\noindent
Hence, the entire process can be described as a linear operator:
\[
T^{(N)} := R F^{(N)} \Lambda^{(N)},
\]
mapping a real non-negative magnitude spectrum \( M^{(N)} \) to the time-reversed complex signal \( \tilde{x}^{(N)} \). An immediate consequence  is seen in the case when the input magnitude spectrum \( M^{(N)} \) is a Kronecker delta function:
\[
M^{(N)}[k] = 
\begin{cases}
1, & \text{if } k = k_0 \text{ for some fixed } k_0, \\
0, & \text{otherwise}.
\end{cases}
\]
In that case, \( X^{(N)} = \Lambda^{(N)} M^{(N)} \) simply extracts one complex exponential term, and \( x^{(N)} \) becomes a complex sinusoid whose magnitude is constant:
\[
|x^{(N)}[n]| = 1 \quad \forall~n.
\]
Consequently, the output magnitude \( |x^{(N)}[n]| \) does not match the input delta-shaped spectrum \( M^{(N)}[n] \), violating the pointwise convergence condition. This illustrates that the conjecture must necessarily exclude sparse spectra or singular inputs like the Kronecker delta and instead be confined to spectrally rich or smooth functions such as those in a dense subset of \( C^\infty([-\pi, \pi]; \mathbb{R}_{\geq 0}) \).

\begin{figure}
    \centering
    \includegraphics[width=3.5in]{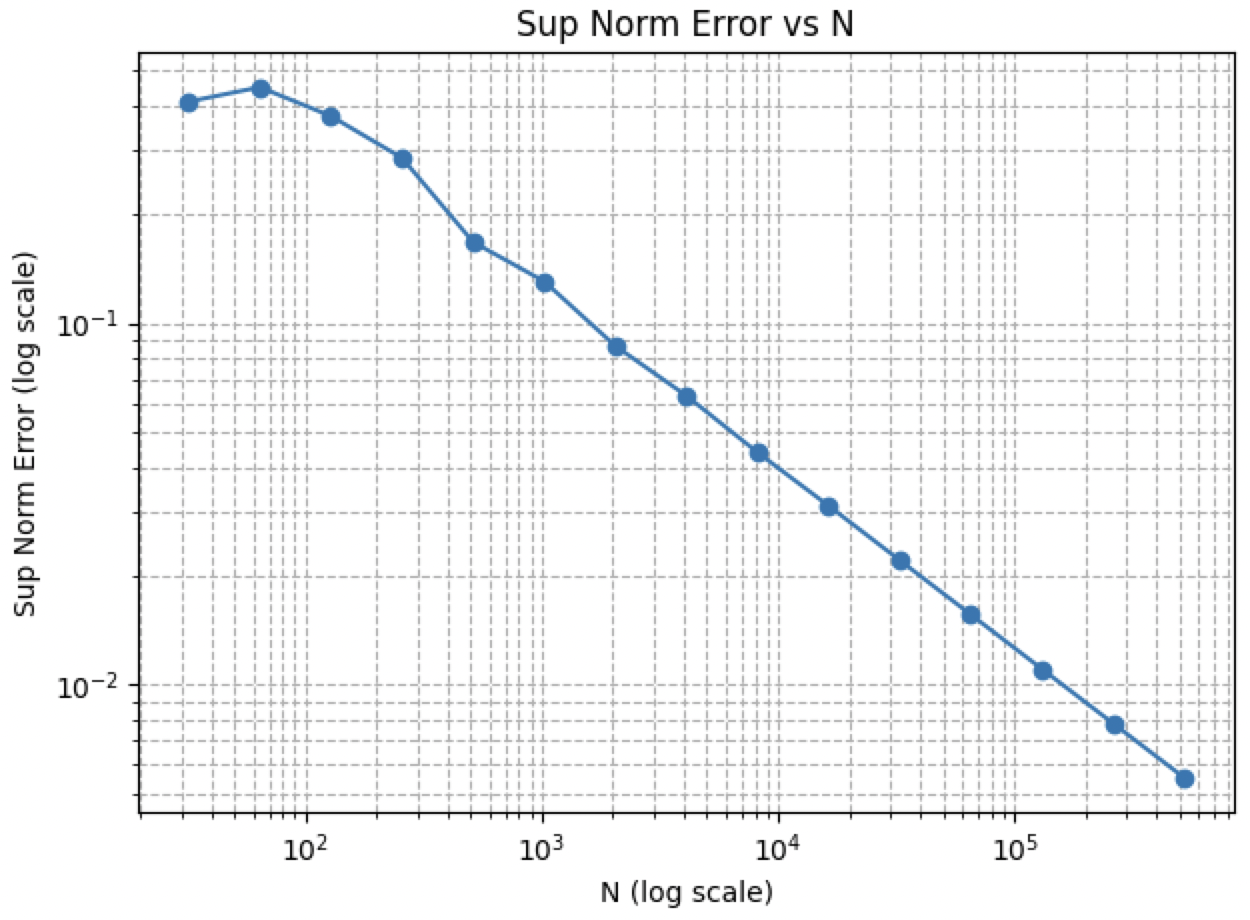}
    \caption{
    Log-log plot of the sup-norm error $\max_k \left|\, \left|x^{(N)}[k]\right| - M^{(N)}[k]\,\right|$
    versus signal length $N$, where $N$ ranges from $2^9$ to $2^{17}$ and includes intermediate values up to $2^{20}$. The observed decay supports the conjecture that the magnitude of the inverse DFT of a spectrum constructed using the Newman phase (with added linear phase $\pi k$) converges uniformly to the input magnitude spectrum as $N \to \infty$.
    }
    \label{fig:diff_vs_n}
\end{figure}

\section*{Numerical Experiments with Random Sum-of-Sinusoids Spectra}

\noindent In order to investigate the empirical behavior of the Newman phase sequence more systematically, a large-scale numerical experiment was performed using a family of randomized input magnitude spectra. Specifically, $1000$ distinct spectra were generated, each constructed as a sum of three sinusoids with random amplitudes, frequencies, and phases. To ensure strict non-negativity, a constant bias term equal to the sum of the absolute values of the sinusoid amplitudes plus an additional constant was added. This guaranteed that the resulting spectra remained positive at all frequency indices.

\medskip

\noindent For each input spectrum, the Newman phase sequence was attached and the inverse discrete Fourier transform (IDFT) was computed to obtain the corresponding time-domain signal. Due to the known effect that the Newman phase sequence induces a time reversal in the signal, the absolute value of the time-domain signal was compared to the flipped version of the original magnitude spectrum. The root mean square (RMS) error between the reconstructed magnitude and the flipped target spectrum was computed for each case.

\medskip

\noindent This procedure was repeated for multiple values of $N$, with $N$ increasing in powers of two from $32$ up to $32 \times 4096$. The resulting RMS errors across all signals were summarized using box plots for each $N$. The observed trend indicated that, despite significant variability in the input signals, the RMS error generally decreased as $N$ increased. However, the rate of convergence was highly non-uniform and remained sensitive to the degree of smoothness and the local structure of the input spectra. These empirical results reinforce the conjecture that the Newman phase sequence works particularly well for a dense subclass of smooth, non-negative functions but that discontinuities or sharp transitions can lead to localized oscillatory artifacts reminiscent of the Gibbs phenomenon.

\bigskip

\section*{Optimization Setup and Observed Behavior}

\noindent To explore whether a better universal phase sequence exists that could achieve uniformly lower reconstruction errors, an explicit optimization problem was formulated. For a fixed $N$, the phase sequence was parameterized as a real vector $\boldsymbol{\theta} \in [0, 2\pi]^N$, so that the applied phase at frequency index $k$ is $\exp(i \theta_k)$. The time-domain signal for each input magnitude spectrum is obtained by attaching the phase, computing the IDFT, and taking the magnitude.

\medskip

\noindent Let $\mathcal{M} = \{ M^{(N)}_1, M^{(N)}_2, \dots, M^{(N)}_{1000} \}$ denote the set of all input magnitude spectra. For each $M^{(N)}_j$, let $x_j^{(N)}(\boldsymbol{\theta})$ denote the resulting time-domain signal, and let $\widetilde{M}^{(N)}_j(\boldsymbol{\theta})$ denote its pointwise magnitude:
\[
\widetilde{M}^{(N)}_j(\boldsymbol{\theta}) = \left| F^{(N)} \cdot \mathrm{diag}( e^{i \boldsymbol{\theta}} ) \cdot M^{(N)}_j \right|, \quad j = 1,\dots,1000.
\]

\medskip

\noindent The optimization problem is then defined as:
\[
\min_{\boldsymbol{\theta} \in [0, 2\pi]^N} \quad \sum_{j=1}^{1000} 
\sqrt{ \frac{1}{N} \sum_{k=0}^{N-1} 
\Big( \widetilde{M}^{(N)}_j[k](\boldsymbol{\theta}) - M^{(N)}_j[N-1-k] \Big)^2 }.
\]

\medskip

\noindent Here, the flipped index $N-1-k$ accounts for the known time-reversal effect of the Newman phase. The bounds $\theta_k \in [0, 2\pi]$ ensure that the phase applied at each frequency lies within the full range of the unit circle.

\medskip

\noindent The solver was implemented using a gradient-based method with the analytically computed Jacobian. In particular, the \texttt{Conjugate Gradient (CG)} method was found to be effective for this smooth, non-convex problem. When initialized with a random phase vector, the optimizer typically reduced the objective by an order of magnitude in a few iterations, but the final value remained significantly higher than that obtained using the Newman phase sequence.

\medskip

\noindent In contrast, when the Newman phase sequence itself was used as the initial guess, the optimizer immediately recognized that it was already at or extremely close to a local minimum, often terminating in zero iterations with no gradient updates needed. This behavior strongly suggests that the Newman phase sequence sits near an isolated global or near-global optimum for this family of smooth spectra. Small random perturbations around the Newman phases resulted in negligible additional improvement, further reinforcing this conclusion.

\medskip

\noindent An immediate exception to the conjectured behavior is the Kronecker delta function: when this trivial magnitude spectrum is used, the operation of attaching the Newman phase, computing the IDFT, and taking the magnitude does not reproduce the same spectrum (even up to time reversal). This highlights that the conjecture is not expected to hold for all possible inputs but rather for a dense subclass of smooth, non-negative spectra.

\medskip

\noindent These numerical results support the conclusion that the Newman quadratic phase construction is not merely convenient but also nearly optimal in minimizing signal-domain oscillations for a wide range of practical input spectra. Furthermore, they illustrate the highly non-convex nature of the phase retrieval landscape, where random initialization rarely produces phases that outperform the explicit closed-form sequence.

\bibliographystyle{IEEEtran}
\bibliography{biblio}
\end{document}